\documentclass[12pt, a4paper]{article}
\usepackage[T1]{fontenc}
\usepackage[french]{babel}
\usepackage{amsmath,amsthm}
\usepackage{amssymb}
\usepackage{amscd}
\usepackage[latin1]{inputenc}
\usepackage{amsfonts}
\usepackage[matrix,arrow,curve]{xy}

\theoremstyle{plain}
\newtheorem{theo}{Th\'eor\`eme}[section]
\newtheorem{lem}[theo]{Lemme}
\newtheorem{prop}[theo]{Proposition}

\theoremstyle{definition}
\newtheorem{df}[theo]{D\'efinition}
\newtheorem{theosd}[theo]{Th\'eor\`eme}
\newtheorem{propsd}[theo]{Proposition}
\newtheorem{rem}[theo]{Remarque}

\newcommand{\Z}{\ensuremath{\mathbb Z/2}}
\newcommand{\FF}{\ensuremath{\mathbb F}}
\newcommand{\KK}{\ensuremath{\mathcal K}}
\newcommand{\Hn}{\ensuremath{H^3_{\mathrm{nr}}}}
\newcommand{\nH}{\ensuremath{H_{\mathrm{nr}}}}
\newcommand{\PF}{\ensuremath{\mathbb P_{\mathbb F_p}}}
\newcommand{\Pk}{\ensuremath{\mathbb P_{k}}}
\newcommand{\QZ}{\ensuremath{\mathbb Q/\mathbb Z}}
\newcommand{\QZl}{\ensuremath{\mathbb Q_l/\mathbb Z_l}}
\newcommand{\emptyline}{$$\,$$}
\newcommand{\eqdef}{\ensuremath{\stackrel{\mathrm{d\acute{e}f}}{=}}}

\title{Sur le groupe de Chow de codimension deux des vari\'et\'es sur les corps finis
}
\author{Alena Pirutka}

\begin{document}

\maketitle
\date{}
\begin{abstract}
En utilisant la construction de Colliot-Th\'el\`ene et Ojanguren, on donne un exemple d'une vari\'et\'e  projective et lisse g\'eom\'etriquement rationnelle $X$,  d\'efinie sur un corps fini $\FF_p$,   telle que d'une part le groupe $\Hn(X, \Z)$ est non nul et, d'autre part, l'application
$CH^2(X)\to CH^2(X\times_{\FF_p}\bar \FF_p)^{\mathrm{Gal}(\bar \FF_p/\FF_p)}$ n'est pas surjective.\\
\begin{center}\textbf{Abstract}\end{center}

Using the construction of Colliot-Th\'el\`ene and Ojanguren, we exhibit an example of a smooth projective geometrically rational variety $X$ defined over a finite field $\FF_p$ with an algebraic closure $\bar \FF_p$, such that the group $\Hn(X, \Z)$ is nonzero and the map
$CH^2(X)\to CH^2(X\times_{\FF_p}\bar \FF_p)^{\mathrm{Gal}(\bar \FF_p/\FF_p)}$ is not surjective.
\end{abstract}
\emptyline

Soit $\FF_p$ un corps fini de cardinal $p$. Soit $\bar \FF_p$ une cl\^oture alg\'ebrique de $\FF_p$ et soit $G=\mathrm{Gal}(\bar \FF_p/\FF_p)$ le groupe de Galois absolu. Soit $X$ une $\FF_p$-vari\'et\'e projective et lisse, g\'eom\'etriquement connexe, de dimension $d$ et soit $\bar X= X\times_{\FF_p} \bar \FF_p$. On dispose d'une application naturelle
$$CH^i(X)\to CH^i({\bar X})^G$$
entre les groupes de Chow des cycles de codimension $i$ sur $X$ (resp. sur $\bar X$) modulo l'\'equivalence rationnelle. Cette application est  surjective pour $i=0,1, d$ (cf. remarque \ref{surj}).
Suivant Geisser \cite{G}, on s'int\'eresse \`a savoir s'il en est ainsi pour $2\leq i<d$.

Dans cet article, on donne un contre-exemple pour $i=2$. Dans ce cas, des arguments de $K$-th\'eorie alg\'ebrique (cf. \cite{K}) permettent de faire un lien entre le conoyau de l'application $CH^2(X)\to CH^2({\bar X})^G$ et le groupe de cohomologie non ramifi\'ee $H^3_{\mathrm{nr}}(X, \QZl(2))$. On montre qu'il suffit d'assurer que ce dernier groupe est non nul (cf. section \ref{kth}). Pour ce faire, les techniques d\'evelopp\'ees par Colliot-Th\'el\`ene et Ojanguren \cite{CTOj} sont disponibles.  En utilisant leur m\'ethode, on construit ainsi (cf. section \ref{cex}) une vari\'et\'e projective lisse $X$ g\'eom\'etriquement connexe d\'efinie sur un corps fini $\FF_p$ convenable, telle que \begin{center}l'application
$CH^2(X)\to CH^2({\bar X})^G$ \textit{n'est pas surjective}.
 \end{center}Plus pr\'ecis\'ement, $X$ est une vari\'et\'e g\'eom\'etriquement rationnelle de dimension $5$, admettant un morphisme vers $\PF^2$ \`a fibre g\'en\'erique une quadrique lisse <<voisine>> de Pfister. Notre m\'ethode permet d'obtenir de tels exemples sur des corps finis $\FF_p$ pour une infinit\'e de nombres premiers $p$.\\

\textbf{Remerciements :} Je tiens \`a exprimer ma profonde gratitude  \`a mon directeur de th\`ese, Jean-Louis Colliot-Th\'el\`ene, pour m'avoir sugg\'er\'e d'utiliser les m\'ethodes de \cite{CTOj} et  m'avoir introduite dans le sujet; sans ses nombreux conseils et r\'eponses cet article n'aurait pas pu voir le jour.

\section{\large{Notations et rappels}}

\subsection{Notations}

\'Etant donn\'e un corps $k$, on note $k^*$ le groupe multiplicatif $k-\{0\}$, $\bar k$ une cl\^oture s\'eparable de $k$ et $G=\mathrm{Gal}(\bar k/k)$ le groupe de Galois absolu. On note $\FF_p$ le corps fini de cardinal $p$.

Si $X$ est une vari\'et\'e alg\'ebrique d\'efinie sur un corps $k$, on note $\bar X=X\times_k\bar k$. Si $X$ est int\`egre, on note $k(X)$ son corps des fractions et si $X$ est g\'eom\'etriquement int\`egre, on note $\bar k(X)$ le corps des fractions de $\bar X$. On dit que $X$ est \textit{$k$-rationnelle} si $X$ est birationnelle \`a $\mathbb P^n_{k}$ et on dit que $X$ est \textit{g\'eom\'etriquement rationnelle} si $\bar X$ est $\bar k$-rationnelle.

Pour une $k$-vari\'et\'e int\`egre $X$ et $i$ un entier, on note $X^{(i)}$ l'ensemble des points de $X$ de codimension $i$ et on note $CH^i(X)$ le groupe des cycles de codimension $i$ modulo l'\'equivalence rationnelle.\\

Si  $A$ est un groupe ab\'elien et $n$ est un entier, on note $A[n]$ le sous-groupe de $A$ form\'e par les \'el\'ements annul\'es par $n$. Pour $l$ un nombre premier, on note $A\{l\}$ le sous-groupe de torsion $l$-primaire.

Pour $M$ un $G$-module continu discret on note $H^i(k, M)=H^i(G, M)$   le $i$-\`eme groupe de cohomologie galoisienne et on note $M^G=H^0(k,M)$ le sous-groupe form\'e par les \'el\'ements invariants par $G$.

\subsection{Rappels de cohomologie \'etale}
\'Etant donn\'es un corps $k$ et un entier $n$ inversible sur $k$, on note $\mu_{n}$ le $k$-sch\'ema en groupes (\'etale) des racines $n$-i\`emes de l'unit\'e. Pour $j$ un entier positif, on note $\mu_{n}^{\otimes j}=\mu_{n}\otimes\ldots\otimes\mu_{n}$ ($j$ fois). On pose $\mu_{n}^{\otimes j}=Hom_{k-gr}(\mu_{n}^{\otimes (-j)}, \mathbb Z/n)$ si $j$ est n\'egatif et $\mu_{n}^{\otimes 0}=\mathbb Z/n$.  Ces $k$-sch\'emas en groupes donnent des faisceaux \'etales, not\'es encore $\mu_{n}^{\otimes j}$, sur toute $k$-vari\'et\'e $X$. On note $H^i(X,\mu_n^{\otimes j})$ les groupes de cohomologie \'etale de $X$ \`a valeurs dans $\mu_n^{\otimes j}$. Lorsque $n=2$, on a un isomorphisme canonique  $\mu_{2}^{\otimes j}\stackrel{\sim}{\to}\Z$ pour tout $j$.

\df{Pour $X$ une $k$-vari\'et\'e int\`egre, un entier naturel $j\geq 1$ et un entier relatif $i$, on d\'efinit les groupes de cohomologie \textit{non ramifi\'ee}
$$\nH^j(X, \mu_n^{\otimes i})\eqdef\nH^j(k(X)/k, \mu_n^{\otimes i})=\bigcap\limits_A\mathrm{Ker}[H^j(k(X), \mu_n^{\otimes i})\stackrel{\partial_{j,A}}{\to}H^{j-1}(k_A, \mu_n^{\otimes i-1})].$$
Dans cette formule, $A$ parcourt les anneaux de valuation discr\`ete de rang un, de  corps des fractions $k(X)$, contenant le corps $k$. Le corps r\'esiduel d'un tel anneau $A$ est not\'e $k_A$ et l'application $\partial_{j,A}$ est l'application r\'esidu. \\}

Lorsque $X$ est propre et lisse, les r\'esultats de Bloch et Ogus permettent d'identifier le groupe $\nH^j(X, \mu_n^{\otimes i})$ au groupe de cohomologie de Zariski $H^0(X, \mathcal H^j(\mu_n^{\otimes i}))$, où $\mathcal H^j(\mu_n^{\otimes i})$ d\'esigne le faisceau de Zariski sur $X$ associ\'e au pr\'efaisceau $U\mapsto H^j(U,\mu_n^{\otimes i})$ (cf. \cite{CT}).

On note $H^i(X, \QZ(j))$, resp. $\nH^{i}(X, \QZ(j))$  (resp. $H^i(X,\mathbb Q_l/\mathbb Z_l(j))$, resp.\\
 $\nH^{i}(X,\mathbb Q_l/\mathbb Z_l(j))$) la limite inductive des groupes $H^i(X,\mu_n^{\otimes j})$, resp. $\nH^{i}(X, \mu_n^{\otimes j})$  lorsque $n$ varie parmi les entiers (resp. parmi les puissances d'un nombre premier $l$, $l\neq \mathrm{car}.k$).\\

 On note $\mathbb G_m$ le groupe multiplicatif sur un sch\'ema $X$ et le faisceau \'etale ainsi d\'efini. On \'ecrit $\mathrm{Br}\,X=H^2_{\acute{e}t}(X,\mathbb G_m)$ pour le groupe de Brauer cohomologique de $X$ et $\mathrm{Pic}(X)=H^1_{\mathrm{Zar}}(X, \mathcal O_X^*)\simeq H^1_{\acute{e}t}(X,\mathbb G_m)$ pour le groupe de Picard.

\subsection{Rappels de $K$-th\'eorie}
Pour $X$ un sch\'ema noeth\'erien et $j$ un entier positif on note $\KK_j$ le faisceau de Zariski associ\'e au pr\'efaisceau $U\mapsto K_j(H^0(U, \mathcal O_U))$, le groupe $K_j(A)$ \'etant celui associ\'e par Quillen \cite{Q} \`a l'anneau $A$.

Lorsque $X$ est une vari\'et\'e lisse sur un corps $k$, la conjecture de Gersten, \'etablie par Quillen \cite{Q}, permet de calculer les groupes de cohomologie de Zariski $H^i(X, \KK_j)$ comme les groupes de cohomologie du complexe de Gersten. Lorsque $j=2$, qui est le cas qui nous int\'eresse dans la suite, ce complexe s'\'ecrit
$$K_2k(X)\stackrel{d_2}{\to}\bigoplus\limits_{x\in X^{(1)}}k(x)^*\stackrel{d_1}{\to} \bigoplus\limits_{x\in X^{(2)}} \mathbb Z,$$
où l'application $d_2$ est donn\'ee par le symbole mod\'er\'e et l'application $d_1$ est obtenue par la somme des fl\`eches diviseurs apr\`es normalisation des vari\'et\'es consid\'er\'ees.
On a ainsi
$H^0(X, \KK_2)=\mathrm{Ker}\,d_2$ et
$H^1(X, \KK_2)=\mathrm{Ker}\,d_1/\mathrm{Im}\,d_2.$

\'Etant donn\'e un corps $k$, le groupe $K_2k$ co\"incide avec le groupe  de $K$-th\'eorie de Milnor $K_2^Mk$, quotient de $k^{*}\otimes_{\mathbb Z} k^*$ par le sous-groupe engendr\'e par les \'el\'ements $a\otimes b$ avec $a+b=1$.

Cette description permet de voir que pour $X$ une vari\'et\'e lisse sur un corps $k$ on a une fl\`eche naturelle
\begin{equation}\label{ph}\mathrm{Pic}(X)\otimes  k^*\to H^1(X, \KK_2).
\end{equation} En effet, on a le diagramme commutatif suivant
$$\begin{CD}  k(X)^*\otimes k^*@>>> \bigoplus\limits_{x\in X^{(1)}}k^*@>>>  \mathrm{Pic}(X)\otimes  k^*@>>> 0\\
@VVV @V\phi VV  @. @. \\
K_2k(X)@>d_2>>\bigoplus\limits_{x\in X^{(1)}}k(x)^*@>d_1>> \bigoplus\limits_{x\in X^{(2)}} \mathbb Z
\end{CD}$$
où la premi\`ere ligne est obtenue \`a partir de la suite exacte  d\'efinissant le groupe $\mathrm{Pic}(X)$ par tensorisation avec $k^*$. On v\'erifie que la compos\'e $d_1\circ \phi$ vaut z\'ero, ce qui permet de d\'efinir la fl\`eche $\mathrm{Pic}(X)\otimes  k^*\to H^1(X, \KK_2)$ par  chasse au diagramme.\\

\section{\large{Comparaison entre groupes de Chow en codimension deux et cohomologie non ramifi\'ee en degr\'e trois}}\label{kth}

Dans cette partie on donne la preuve du th\'eor\`eme suivant :

\theosd\label{comp}{\textit{Soit $X$ une $\mathbb Q$-vari\'et\'e projective et lisse, g\'eom\'etriquement rationnelle. Pour presque tout nombre premier $p$, la r\'eduction $X_p$ de $X$ modulo $p$ est bien d\'efinie et est une $\mathbb F_p$-vari\'et\'e projective et lisse, g\'eom\'etriquement rationnelle, telle que $$\Hn(X_p, \QZl(2))\stackrel{\simeq}{\to}\mathrm{Coker}[CH^2(X_p)\to CH^2(\bar X_p)^G]\{l\}$$
pour tout nombre premier $l$, $(l,p)=1$. \\}}

Pour d\'emontrer ce th\'eor\`eme, on utilise le r\'esultat suivant de B. Kahn (\cite{K}, Th.1 et corollaire p.$3$) :

\theosd\label{Kahn}(\cite{K}){\textit{ Soit $k$ un corps de caract\'eristique $p\geq 0$, de dimension cohomologique au plus $3$. Soit $X$ une $k$-vari\'et\'e projective et lisse. Supposons que les conditions suivantes sont satisfaites :
\begin{itemize}
\item [(i)] $K_2\bar k\stackrel{\sim}{\to} H^0(\bar X, \KK_2)$;
\item [(ii)] le groupe $\nH^3(\bar X,\QZ(2))$ est nul, resp. de torsion $p$-primaire si $\mathrm{car}.k>0$.\\
\end{itemize}
Alors on a une suite exacte naturelle, resp. exacte \`a la $p$-torsion pr\`es si $\mathrm{car}.k>0$
\begin{multline}\label{sf}H^1(k,H^1(\bar X,\KK_2))\to\mathrm{Coker}[H^3(k, \QZ(2))\to \Hn(X, \QZ(2))]\to\\\to\mathrm{Coker}[CH^2(X)\to CH^2(\bar X)^G]\to H^2(k,H^1(\bar X,\KK_2)).\end{multline}}}

\rem{Les groupes de cohomologie \`a coefficients dans $\QZ(2)$ sont bien d\'efinis en caract\'eristique positive (cf.\cite{K}). \\}

Il est ainsi n\'ecessaire  de v\'erifier les hypoth\`eses $(i)$ et $(ii)$ pour une vari\'et\'e  g\'eom\'etriquement rationnelle $X$. Les \'enonc\'es suivants, cas particuliers de \cite{CT}, 2.1.9 (cf. aussi 4.1.5), sont bien connus.

\propsd\label{is}{\textit{Soit $k$ un corps. Soit $X$ une $k$-vari\'et\'e projective et lisse, $k$-rationnelle. Alors
\begin{itemize}
\item [(i)] $\nH^j(X, \mu_n^{\otimes i})\simeq H^j(k, \mu_n^{\otimes i})$ pour tout $j\geq 1$;
\item[(ii)] l'application naturelle $K_2k\to H^0(X, \KK_2)$ est un isomorphisme;
\item[(iii)] le groupe $\mathrm{Pic}(X)$ est libre de type fini.\\
\end{itemize}
}}

L'\'enonc\'e suivant permet de comprendre le module galoisien $H^1(\bar X, \KK_2)$. 

\prop\label{iso0}{Soit  $k$ un corps de caract\'eristique nulle. Le noyau $K(X)$ et le conoyau $C(X)$ de l'application $\mathrm{Pic}(X)\otimes k^*\to H^1(X, \KK_2)$ sont des invariants birationnels des $k$-vari\'et\'es projectives et lisses, g\'eom\'etriquement int\`egres. En particulier, l'application $\mathrm{Pic}(X)\otimes k^*\to H^1(X, \KK_2)$  est un isomorphisme pour $X$  une vari\'et\'e projective et lisse, $k$-rationnelle. }

\proof{Consid\'erons le complexe de groupes ab\'eliens

$$\mathrm{Pic}(X) \otimes k^* \to H^1(X,K_2).$$

Ce complexe est fonctoriel contravariant pour les morphismes dominants
de vari\'et\'es projectives et lisses.
Le noyau $K(X)$ et le conoyau $C(X)$ de $\mathrm{Pic}(X) \otimes k^* \to H^1(X,K_2)$
sont alors des foncteurs contravariants pour de tels morphismes. Soit $F$ l'un de ces foncteurs.

Soient $X,Y$ deux vari\'et\'es projectives et lisses, g\'eom\'etriquement int\`egres. Montrons  qu'un morphisme birationnel $X\to Y$  induit un isomorphisme $F(Y)\to F(X)$. D'apr\`es Hironaka, il existe deux $k$-vari\'et\'es projectives et lisses $X'$ et $Y'$ et un diagramme commutatif

$$\xymatrix{
X'\ar[r]\ar[d] &Y'\ar[d]\ar[ld] &\\
X\ar[r] & Y &
}
$$
où les fl\`eches verticales sont  des suites d'\'eclatements de centres lisses. D'apr\`es le lemme ci-dessous, $F(X)$ est isomorphe à $F(X')$, respectivement $F(Y)$ est isomorphe à $F(Y')$. On en d\'eduit par fonctorialit\'e que $F(X)$ est isomorphe à $F(Y)$.

Si maintenant on a une application rationnelle $X\dashrightarrow Y$, on utilise Hironaka pour trouver une vari\'et\'e projective et lisse $Z$ avec $Z\to X$ et $Z\to Y$ deux morphismes birationnels. D'apr\`es ce qui pr\'ec\`ede, $F(X)\simeq F(Z)\simeq F(Y)$. Ainsi $F(X)$ est un invariant birationnel des $k$-vari\'et\'es projectives et lisses, g\'eom\'etriquement int\`egres.

Le fait que  $F(\mathbb P^n_{k})=0$  est bien connu. On \'etablit d'abord que $H^1(\mathbb A^1_k, \KK_2)=0$ et ensuite que $H^1(\mathbb A^n_k, \KK_2)=0$ par des fibrations successives à fibres $\mathbb A^1$. L'\'enonc\'e pour $\mathbb P^n_{k}$ s'en suit par r\'ecurrence, en se restreignant à l'hyperplan à l'infini.  \qed\\}

\rem{On peut montrer plus g\'en\'eralement que pour  $X$ lisse sur un corps, $H^i(\mathbb A^n_X,\mathcal K_j)$ est isomorphe à $H^i( X,\mathcal K_j)$, et donner une expression explicite de $H^i(\mathbb P^n_X,\mathcal K_j)$ en termes de $K$-cohomologie de $X$, cf. \cite{Sh}.\\}

\lem\label{lemec}{Soit $k$ un corps. Soit $X$ une $k$-vari\'et\'e projective et lisse, g\'eom\'etriquement int\`egre. Soit $Z\subset X$ une sous-vari\'et\'e int\`egre, projective et lisse, de codimension au moins $2$ et soit $\pi:X'\to X$ l'\'eclatement de $X$ le long de $Z$. Alors les applications $K(X)\to K(X')$, respectivement  $C(X)\to C(X')$, sont des isomorphismes.}
\proof{  Soit $Z'$ le diviseur exceptionnel de $X'$ et soit $U=X\setminus Z\simeq X'\setminus Z'.$ Supposons d'abord que $Z$ est de codimension $2$. On a les suites exactes horizontales de complexes verticaux :

\small{
$$\begin{CD}
@. 0@>>>  K_2k(X) @>>>K_2k(U) @>>>0\\
@. @VVV   @VVV  @VVV   @.\\
@. 0@>>> \bigoplus\limits_{x\in X^{(1)}} k(x)^*@>>> \bigoplus\limits_{x\in U^{(1)}}k(x)^*@>>>0\\
@. @VVV   @VVV  @VVV   @.\\
0 @>>> \mathbb Z @>>>\bigoplus\limits_{x\in X^{(2)}}\mathbb Z@>>>\bigoplus\limits_{x\in U^{(2)}}\mathbb Z@>>>0
\end{CD}$$}

\normalsize et

\small{
$$\begin{CD}
@. 0@>>>  K_2k(X') @>>>K_2k(U) @>>>0\\
@. @VVV   @VVV  @VVV   @.\\
0@>>> k(Z')^*@>>> \bigoplus\limits_{x\in {X'}^{(1)}} k(x)^*@>>> \bigoplus\limits_{x\in U^{(1)}}k(x)^*@>>>0\\
@. @VVV   @VVV  @VVV   @.\\
0 @>>> \bigoplus\limits_{x\in {Z'}^{(1)}}\mathbb Z @>>>\bigoplus\limits_{x\in {X'}^{(2)}}\mathbb Z@>>>\bigoplus\limits_{x\in U^{(2)}}\mathbb Z@>>>0.
\end{CD}$$}
\normalsize

On a ainsi des suites longues induites en cohomologie  :
\begin{multline}\label{st1}
 0\to H^0(X, \KK_2)\to H^0(U, \KK_2)\to 0\to H^1(X, \KK_2)\to H^1(U, \KK_2)\to  \mathbb Z\to \\\to CH^2(X)\to CH^2(U)\to 0,
\end{multline}

\begin{multline}\label{st2}
 0\to H^0(X', \KK_2)\to H^0(U, \KK_2)\to k^*\to H^1(X', \KK_2)\to H^1(U, \KK_2)\to  \mathrm{Pic}(Z')\to\\\to CH^2(X')\to CH^2(U)\to 0.
\end{multline}

Dans la suite (\ref{st1}), la fl\`eche  $\mathbb Z\to CH^2(X)$ est donn\'ee par $1\mapsto[Z]$. En prenant l'intersection avec un hyperplan g\'en\'eral, on voit que cette fl\`eche est injective. Ainsi l'application $H^1(X, \KK_2)\to H^1(U, \KK_2)$ est un isomorphisme. 

Si $Z$ est de codimension plus grande que $2$, on a encore la suite (\ref{st2}) et les groupes $H^1(X, \KK_2)$ et $H^1(U, \KK_2)$ sont isomorphes, car ils ne dependent que de points de codimension au plus $2$. 

Par fonctorialit\'e, on a  le diagramme commutatif suivant  :
$$\xymatrix{ H^0(X, \KK_2)\ar[r]^{\simeq}\ar[d]&H^0(U, \KK_2)\ar@{=}[d]&\\
H^0(X', \KK_2)\ar@{^{(}->}[r]&H^0(U, \KK_2).&}$$
Ainsi l'application $H^0(X', \KK_2)\to H^0(U, \KK_2)$ est un isomorphisme. En utilisant la suite (\ref{st2}), on obtient  le diagramme commutatif de suites exactes suivant :
$$\begin{CD} @. @. H^1(X, \KK_2)@>\simeq>>H^1(U, \KK_2)\\
@. @. @VVV @|  \\
0@>>>k^*@>>>H^1(X', \KK_2)@>>>H^1(U, \KK_2).
\end{CD}$$

On a donc une suite exacte scind\'ee :
$$0\to k^*\to H^1(X', \KK_2)\to H^1(X, \KK_2)\to 0.$$
Ainsi $H^1(X', \KK_2)\simeq H^1(X, \KK_2)\oplus k^*$.
Puisque $\mathrm{Pic}(X')=\mathrm{Pic}(X)\oplus \mathbb Z\cdot[Z']$, on en d\'eduit l'\'enonc\'e du lemme.
\qed\\}

\paragraph{Preuve du th\'eor\`eme \ref{comp}.}
Puisque $X$ est une $\mathbb Q$-vari\'et\'e g\'eom\'etriquement rationnelle, il existe une extension finie $K/\mathbb Q$ et  une $K$-vari\'et\'e projective et lisse $Z$, deux morphismes $K$-birationnels $Z\to X_K$ et $Z\to \mathbb P^n_K$, et deux diagrammes commutatifs  $$\xymatrix{
Z'\ar[r]\ar[d] &X'\ar[d]\ar[ld] &\\
Z\ar[r] & X_K &
}
\mbox{ et }\;\xymatrix{
Z''\ar[r]\ar[d] &Y\ar[d]\ar[ld] &\\
Z\ar[r] &\mathbb P^{n}_{K} &
}
$$
où les fl\`eches verticales sont  des suites d'\'eclatements de centres lisses.  Pour presque toute place $v$ de $K$, les centres d'\'eclatements admettent des r\'eductions lisses et on a des diagrammes commutatifs induits sur le corps r\'esiduel $k(v)$, donc sur  $\bar \FF_p$, $p=\mathrm{car.}k(v)$.

L'argument donn\'e dans la preuve de la proposition \ref{iso0} s'applique alors : pour presque tout nombre premier $p$, la r\'eduction $X_p$ de $X$ modulo $p$ est bien d\'efinie et est une $\mathbb F_p$-vari\'et\'e projective et lisse, g\'eom\'etriquement rationnelle, telle que l'application $\mathrm{Pic}(\bar X_p)\otimes \bar \FF_p^*\to H^1(\bar X_p, \KK_2)$ est un isomorphisme.

Montrons que les hypoth\`eses du th\'eor\`eme \ref{Kahn} sont satisfaites pour une telle r\'eduction $X_p$. D'apr\`es la proposition \ref{is},  $K_2\bar \FF_p\stackrel{\sim}{\to} H^0(\bar X_p, \KK_2)$ car $X_p$ est g\'eom\'etriquement rationnelle. De m\^eme, $\nH^3(\bar X_p,\QZl(2))= H^3(\bar \FF_p,\QZl(2))=0$ car $\bar \FF_p$ est s\'eparablement clos. 

Montrons ensuite que le groupe $H^i(\mathbb F_p,H^1(\bar X_p,\KK_2))\simeq H^i(\mathbb F_p, \mathrm{Pic}(\bar X_p)\otimes \bar {\FF_p}^*)$ est nul pour tout $i\geq 1$. D'apr\`es la proposition \ref{is},  le $\mathbb Z$-module $\mathrm{Pic}(\bar X_p)$ est libre de type fini.   Consid\'erons une extension finie galoisienne $L/\mathbb F_p$ qui  d\'eploie $\mathrm{Pic}(\bar X_p)$. Consid\'erons la suite de restriction-inflation :
$$0\to H^1(\mathrm{Gal}(L/\mathbb F_p), \mathrm{Pic}X_{p,L}\otimes L^*)\to H^1(\mathbb F_p, \mathrm{Pic}\bar X_p\otimes \mathbb F_p^*)\to H^1(\mathrm{Gal}(\bar{\mathbb F_p}/L), \mathrm{Pic}\bar X_p\otimes \bar {\mathbb F_p}^*).$$ On a $H^1(\mathrm{Gal}(\bar {\mathbb F_p}/L), \mathrm{Pic}\bar X_p\otimes \bar {\mathbb F_p}^*)=0$ d'apr\`es le th\'eor\`eme de Hilbert $90$. Puisque la dimension cohomologique de $\mathbb F_p$ est  $1$,  $H^1(\mathrm{Gal}(L/\mathbb F_p), \mathrm{Pic}X_{p,L}\otimes L^*)=0$ (cf. \cite{Se}, p.170). On a donc $H^i(\mathbb F_p,\mathrm{Pic}(\bar X_p)\otimes \bar {\mathbb F_p}^*)=0$ pour tout $i\geq 1$.
 
Notons que $H^3(\FF_p, \QZl(2))=0$ car la dimension cohomologique de $\FF_p$ est  $1$. La suite (\ref{sf}) donne alors un isomorphisme :$$\Hn(X_p, \QZl(2))\stackrel{\simeq}{\to}\mathrm{Coker}[CH^2(X_p)\to CH^2(\bar X_p)^G]\{l\}.$$ \qed

\rem{En utilisant des r\'esultats de \cite{CTR}, on peut montrer un \'enonc\'e plus g\'en\'eral (cf. aussi \cite{CTV}) :

\textit{Soit $k$ un corps de caract\'eristique $p\geq 0$. Soit $l$ un nombre premier, $(l,p)=1$. Soit $X$ une $k$-vari\'et\'e projective et lisse, g\'eom\'etriquement rationnelle. Supposons que la dimension cohomologique de $k$ est au plus $1$. Alors on a un isomorphisme 
$\Hn(X, \QZl(2))\stackrel{\simeq}{\to}\mathrm{Coker}[CH^2(X)\to CH^2(\bar X)^G]\{l\}.$}

Indiquons comment on proc\`ede pour la preuve. On utilise d'abord \cite{CTR}, 2.12 et 2.14 pour montrer que pour $X$ une $k$-vari\'et\'e projective et lisse, g\'eom\'etriquement rationnelle, le noyau $K(X)$ et le conoyau $C(X)$ de l'application  $\mathrm{Pic}(\bar X)\otimes \bar k^*\to H^1(\bar X, \KK_2)$ sont uniquement divisibles par tout entier premier \`a $p=\mathrm{car}k$. Ici on a encore $H^i(k,\mathrm{Pic}(\bar X)\otimes \bar k^*)=0$ pour  $i\geq 1$, et on conclut comme dans la preuve du th\'eor\`eme  \ref{comp}.

L'utilisation de r\'esultats de \cite{CTR}  fait appel à des techniques tr\`es \'elabor\'ees. Pour ce dont on a besoin dans la suite, le th\'eor\`eme \ref{comp} suffit.  }

\section{\large{L'exemple}}\label{cex}

Dans cette partie, pour une infinit\'e de nombres premiers $p$, on construit une vari\'et\'e projective et lisse g\'eom\'etriquement rationnelle $X$, d\'efinie sur le
corps fini $\mathbb F_p$,  telle que l'application
$$CH^2(X)\to CH^2(\bar X)^{G}$$
n'est pas surjective.

\rem\label{surj}{Si $X$ est  une vari\'et\'e projective et lisse, g\'eom\'etriquement int\`egre, d\'efinie sur le
corps fini $\mathbb F_p$, l'application $CH^i(X)\to CH^i(\bar X)^{G}$ est surjective pour $i=0,1,d$. Le cas $i=0$ est imm\'ediat. Pour $i=1$,  $\mathrm{Pic}(X)\stackrel{\sim}{\to}CH^1(X)$ car $X$ est lisse. Puisque $X$ est projective et g\'eom\'etriquement int\`egre, $\bar k[X]^*=\bar k^*$ et la suite spectrale $E_2^{pq}=H^p(G,H^q(\bar X,\mathbb G_m))\Rightarrow H^{p+q}(X,\mathbb G_m)$ donne une suite exacte $$0\to \mathrm{Pic}(X)\to  \mathrm{Pic}(\bar X)^G\to H^2(G,\bar k^*)\to \mathrm{Br}\,X.$$ Puisque le groupe $H^2(G,\bar k^*)=\mathrm{Br}\,k$ est nul pour un corps fini, on a la surjectivit\'e pour $i=1$. Plus g\'en\'eralement, il en est ainsi  pour toute vari\'et\'e $X$ projective et lisse, g\'eom\'etriquement int\`egre, avec un point rationnel, d\'efinie sur un corps $k$ quelconque : pour une telle vari\'et\'e l'application $H^2(G,\bar k^*)\to \mathrm{Br}\,X$ est injective.

Pour $i=d$, c'est-\`a-dire dans le cas de z\'ero-cycles, on sait que $X$ poss\`ede un z\'ero-cycle de degr\'e $1$ d'apr\`es les estimations de Lang-Weil (cf.\cite{LW}). Il suffit donc de voir que l'application entre les groupes de Chow de z\'ero-cycles de degr\'e z\'ero $A_0(X)\to A_0(\bar X)^G$ est surjective. Ceci r\'esulte de la comparaison de ces derniers groupes avec les points rationnels (resp. les $\bar \FF_p$-points) de la vari\'et\'e d'Albanese $\mathrm{Alb}_X$ de $X$.  En effet, l'application  $A_0(X)\to \mathrm{Alb}_X(\FF_p)$ est surjective (cf. \cite{KS}, Prop. 9, p.274), et l'application $A_0(\bar X)\to \mathrm{Alb}_X(\bar \FF_p)$ est un isomorphisme (cf.\cite{R} et \cite{M}).\\}

D'apr\`es le th\'eor\`eme \ref{comp}, si $X$ est g\'eom\'etriquement rationnelle, il
suffit d'assurer que le groupe $\Hn(X, \QZl(2))$ est non nul pour un certain
nombre premier $l$, $l\neq p$.
Dans l'article \cite{CTOj}, Colliot-Th\'el\`ene et Ojanguren construisent de tels
exemples sur le corps des complexes pour $l=2$.
Les vari\'et\'es qu'ils construisent sont unirationnelles (c'est-\`a-dire,
domin\'ees par un ouvert de l'espace projectif). Via la proposition
\ref{is}, ils obtiennent ainsi des exemples de vari\'et\'es unirationnelles non
rationnelles. Dans la suite, on utilise la m\'ethode de \cite{CTOj} pour
produire des exemples sur les corps finis.\\

La strat\'egie est la suivante :
\begin{enumerate}
\item On consid\`ere une quadrique projective et lisse $Q$ sur le corps
$F=\FF_p(x,y)$, $p\neq 2$, d\'efinie dans $\mathbb P^4_F$ par une \'equation homog\`ene
\begin{equation}\label{eqq}x_0^2-ax_1^2-fx_2^2+afx_3^2-g_1g_2x_4^2=0
\end{equation} où $a\in \FF_p$ est une constante et $f,g_1,g_2\in F$. La quadrique $Q$ admet un point rationnel sur $\bar \FF_p(x,y)$, elle est donc $\bar \FF_p(x,y)$-rationnelle.
\item On donne des conditions n\'ecessaires sur les coefficients dans (\ref{eqq}) pour que le
cup-produit $(a, f, g_1)$ soit non nul dans $\Hn(Q, \Z)$.
\item On v\'erifie que l'on peut trouver $a\in \mathbb Z$ et $f,g_1,g_2\in\mathbb Z(x,y)$ tels que leurs r\'eductions modulo $p$ v\'erifient les conditions de l'\'etape pr\'ec\'edente pour le corps $\FF_p$ pour une infinit\'e de nombres premiers $p$. Par Hironaka, on  trouve une vari\'et\'e  projective et lisse $X$ d\'efinie sur $\mathbb Q$, admettant une fibration sur $\mathbb P^2_{\mathbb Q}$ de fibre g\'en\'erique  la quadrique d\'efinie par (\ref{eqq}). Pour presque tout $p$, $X$ admet une r\'eduction  $X_p$  modulo $p$ qui  est lisse sur $\FF_p$ et pour une infinit\'e de premiers $p$ le groupe $\Hn(X_p, \Z)$ est ainsi non nul.
\end{enumerate}

\subsection{Cohomologie des quadriques}
On commence par citer un r\'esultat d'Arason \cite{A} sur la cohomologie des quadriques.

Soit $k$ un corps, $\mathrm{car.}k\neq 2.$ Soit $\phi$ une forme quadratique non d\'eg\'en\'er\'ee de dimension $m$ d\'efinie sur $k$. On note $X_{\phi}$ la quadrique projective et lisse dans $\mathbb P^{m-1}_k$ d\'efinie par $\phi$. On appelle \textit{$n$-forme de Pfister} sur $k$ une forme quadratique de type $\langle 1, -a_1\rangle\otimes\ldots\otimes\langle 1,-a_n\rangle$, $a_i\in k^*$. Une forme quadratique non d\'eg\'en\'er\'ee $\phi$ est dite \textit{<<voisine de Pfister>>} s'il existe une forme de Pfister $\phi'$ sur $k$ et $a\in k^*$ tels que $\phi$ soit une sous-forme de $a\phi'$ et que la dimension de $\phi$ soit strictement sup\'erieure \`a la moiti\'e de la dimension de $\phi'$.

\theosd\label{ta}(cf. \cite{A}){ \textit{Soit $k$ un corps, $\mathrm{car.}k\neq 2.$ Soit $\phi$ une forme quadratique d\'efinie sur $k$, voisine d'une $3$-forme de Pfister $\langle 1, -a_1\rangle\otimes\langle 1,-a_2\rangle\otimes\langle 1,-a_3\rangle$. Alors
\begin{equation}\label{eta}
 \mathrm{ker}[H^3(k,\Z)\to H^3(k(X_{\phi}), \Z)]=\Z(a_1,a_2, a_3),
\end{equation}
chaque $a_i$ \'etant identifi\'e \`a sa classe dans $H^1(k,\Z)\simeq k^*/k^{*2}$.\\}}

Soit $Q$ la quadrique d\'efinie sur le corps $F=\FF_p(x,y)$, $p\neq 2$, par l'\'equation homog\`ene (\ref{eqq}). D'apr\`es le th\'eor\`eme d'Arason, $$\mathrm{ker}[H^3(F,\Z)\to H^3(F(Q), \Z)]=\Z(a,f, g_1g_2).$$ Pour trouver un \'el\'ement non nul dans $\Hn(Q,\Z)$, on peut ainsi essayer de chercher un \'el\'ement de $H^3(F,\Z)$, diff\'erent de $(a,f, g_1g_2)$ et qui devient non ramifi\'e dans $H^3(F(Q), \Z)$. On va choisir les \'el\'ements $a,f,g_1\mbox{ et } g_2$ pour que l'\'el\'ement $(a, f, g_1)$ convienne.

Faisons d'abord quelques rappels sur les calculs de r\'esidus.

\propsd\label{prres}(\cite{CTOj}, 1.3 et 1.4){\textit{ Soit $A$ un anneau de valuation discr\`ete, de corps des fractions $K$ et de corps r\'esiduel $k$. Soit $j\geq 1$ un entier. \begin{enumerate}
\item Soit $\alpha\in H^j(A, \Z)$ et soit $\alpha_0\in H^j(k,\Z)$ son image par l'application de r\'eduction. Soit $b\in K^*$ de valuation $m$ dans $A$ et soit $\beta$ la classe de $b$ dans $H^1(K,\Z)$. Alors $\partial_A(\alpha\cup\beta)=m\alpha_0$.
\item Soit $\alpha\in H^j(K, \Z)$ et soit $b\in A^*$ dont la classe est un carr\'e dans $k$. Soit $\beta$ la classe de $b$ dans $H^1(K,\Z)$. Alors $\partial_A(\alpha\cup\beta)=0$.\\
\end{enumerate}}}

On d\'ecrit ensuite les conditions qu'on va imposer sur les coefficients de la quadrique $Q$ :

\prop\label{cond}{Soit $F=\FF_p(x,y)$ le corps des fractions rationnelles \`a deux variables sur le corps fini $\FF_p$. Soit  $a\in \FF_p^*\setminus \FF_p^{*2}$ et soient $f,g_1,g_2\in F$ non nuls. Soit $Q$ la quadrique lisse dans $\mathbb P^4_F$ d'\'equation homog\`ene
$$x_0^2-ax_1^2-fx_2^2+afx_3^2-g_1g_2x_4^2=0.$$
Supposons
\begin{enumerate}
\item pour tout $i=1,2$, il existe un anneau de valuation discr\`ete $B_i$ de corps des fractions $F$, tel que $\partial_{B_i}(a, f, g_i)\neq 0$;
\item pour tout anneau de valuation discr\`ete $B$ de  corps des fractions $F$, associ\'e \`a un point de codimension $1$ de $\PF^2$, soit $\partial_B(a, f, g_1)=0$, soit $\partial_B(a, f, g_2)=0$.
\item pour tout anneau de valuation discr\`ete $B$ de corps des fractions $F$, centr\'e en un point ferm\'e $M$ de $\PF^2$, quitte \`a la multiplier par un carr\'e dans $F^*$, l'une au moins des fonctions $f,g_1,g_2$ est inversible en $M$.
\end{enumerate}
Alors l'image $\xi_{F(Q)}$ du cup-produit $\xi=(a, f, g_1)$ dans  $H^3(F(Q),\Z)$ est un \'el\'ement non nul de $H^3_{\mathrm{nr}}(Q,\Z)$. }
\proof{Notons d'abord que  $(a, f, g_1)$ est non nul dans $H^3(F(Q),\Z)$. Sinon, d'apr\`es le th\'eor\`eme \ref{ta}, on a soit $(a, f, g_1)=0$, soit $(a, f, g_1)=(a,f,g_1g_2)$. Ainsi soit $(a, f, g_1)=0$, soit $(a, f, g_2)=0$, contradiction avec la condition 1.

Montrons que
\begin{multline*}\label{hyp}\mbox{pour tout anneau de valuation discr\`ete }B\mbox{ de }F,\\\mbox{ soit }\partial_B(a, f, g_1)=0,\mbox{ soit }\partial_B(a, f, g_2)=0. \quad(*)
 \end{multline*}
Pour un tel anneau $B$ on dispose  d'un morphisme Spec$\,B\to \PF^2$ et les cas 2 et 3 correspondent \`a deux possibilit\'es pour l'image du point ferm\'e de $B$.  La condition $2$ assure  $(*)$ si cette image est un point de codimention $1$ de $\PF^2$.  Sinon l'image du point  Spec$\,k_B$ est un point ferm\'e $M$ de $\PF^2$.  Soit $\mathcal O_M$ l'anneau local de $M$, son corps des fractions est $F$. On dispose d'un morphisme d'anneaux $\mathcal O_M\to B$.

On peut supposer, sans perte de g\'en\'eralit\'e, que la fonction $g_1$ est inversible dans $\mathcal O_M$, quitte \`a la multiplier par un carr\'e.
Ainsi la fonction $g_1$ est inversible dans $B$. Soit $m$ la valuation de $f$ dans $B$. D'apr\`es \ref{prres}.1,
$\partial_B(a, f, g_1)=\partial_B(g_1, a, f)=m(\bar g_1, \bar a)$, où l'on note $\bar g_1$ (resp. $\bar a$)  la classe de $g_1$ (resp. $a$) dans $H^{1}(k_B,\Z)$. Comme $g_1$ et $a$ sont inversibles dans $\mathcal O_M$, ces derni\`eres classes proviennent de classes dans  $H^{1}(k_M,\Z)$. Ainsi
  $(\bar g_1, \bar a)$ provient d'un \'el\'ement de $H^{2}(k_M,\Z)$. Ce dernier groupe est nul, car $k_M$, \'etant un corps fini, est de dimension cohomologique $1$. Ainsi $\partial_B(a, f, g_1)=0$.

Montrons maintenent  que $\xi_{F(Q)}$ est non ramifi\'e. Soit $A$ un anneau de valuation discr\`ete de $F(Q)$ de corps r\'esiduel $k_A$. Si $A$ contient $F$, alors $\xi_{F(Q)}$ provient d'un \'el\'ement de $H^3(A,\Z)$ et son r\'esidu est donc nul. Supposons que $A$ ne contient pas $F$. Alors $B=A\cap F$ est un anneau de valuation discr\`ete de $F$. Soit $k_B$ son corps r\'esiduel. On a le diagramme commutatif suivant (cf. \cite{CTOj}, §1) :
$$\xymatrix{
H^{3}(F(Q),\Z)\ar[r]^{\partial_A} &H^{2}(k_A,\Z) &\\
H^{3}(F,\Z)\ar[r]^{\partial_B}\ar[u]^{\mathrm{res}_{F/F(Q)}} &H^{2}(k_B,\Z).\ar[u]_{e_{B/A}\mathrm{res}_{k_B/k_A}} &
}
$$
D'apr\`es ce qui pr\'ec\`ede, $\partial_B(a, f, g_i)=0$ pour $i=1$ ou pour $i=2$.
Si  $\partial_B(a, f, g_1)=0$, alors $\partial_A(a, f, g_1)=0$  d'apr\`es le diagramme ci-dessus. Supposons que $\partial_B(a, f, g_2)=0$. Ainsi $\partial_A(a, f, g_2)=0$. Comme $(a,f,g_1g_2)$ est nul dans $H^{3}(F(Q),\Z)$, son r\'esidu l'est aussi dans  $H^{2}(k_A,\Z)$. On a donc $\partial_A(a, f, g_1)=\partial_A(a, f, g_1g_2)-\partial_A(a, f, g_2)=0$. Ainsi  $\xi_{F(Q)}$ est non ramifi\'e.  \qed\\}

\subsection{Construction explicite}

On proc\`ede maintenant \`a la construction des exemples.

Soit $k$ un corps. Dans la suite, on va prendre $k=\FF_p$ ou $k=\mathbb Q$. On fixe $x,y,z$ des coordonn\'ees homog\`enes pour $\Pk^2$. Soit $a\in k^*\setminus k^{*2}$. Soient $b_i, c_i, d_i\in k^{*}\setminus\{-1\}$, $i=1,2$, et soit $l_i=b_ix+c_iy+d_iz$. Soient $h_j$, $j=1,\ldots, 8$, les formes lin\'eaires $e_xx+e_yy+e_zz$, $e_x,e_y, e_z\in\{0, 1\}$.\\

On choisit $b_i, c_i, d_i$ de sorte que :
\begin{itemize}
 \item[(i)] Les droites dans $\Pk^2$ donn\'ees par les \'equations $x=0$, $y=0$, $z=0$,  $l_i+h_j=0$, $i=1,2$, $j=1,\ldots, 8$, soient deux \`a deux distinctes.
 \item[(ii)] Pour tous $1\leq j, j'\leq 8$ les trois droites $x=0$, $l_1+h_j=0$, $l_2+h_{j'}=0$ dans $\Pk^2$ sont d'intersection vide.
\item[(iii)] Pour tous $1\leq j, j'\leq 8$ les trois droites $y=0$, $l_1+h_j=0$, $l_2+h_{j'}=0$  dans $\Pk^2$ sont d'intersection vide. \\
\end{itemize}

On prend pour $f,g_1,g_2\in k(\Pk^2)$ les \'el\'ements suivants :
\begin{equation} \label{coef}f=\frac{x}{y},\quad g_1=\frac{\prod_j(l_1+h_j)}{y^8},\quad
g_2=\frac{\prod_j(l_2+h_j)}{z^8}.\end{equation}
\medskip

\rem{Soit $h_j=e_xx+e_yy+e_zz$. Les droites $x=0$, $l_1+h_j=0$ s'intersectent en un seul point $[0: (c_1+e_y) : (d_1+e_z)]$.
 Ainsi les conditions (ii) et (iii) ci-dessus sont \'equivalentes aux conditions suivantes :
\begin{itemize}
\item[(ii')] les ensembles $$\{[(c_1+e_y) : (d_1+e_z)], e_y, e_z\in \{0,1\}\}\mbox{ et }\{[(c_2+e_y) : (d_2+e_z)], e_y, e_z\in \{0,1\}\}$$ sont d'intersection vide;
\item[(iii')]  de m\^eme,  $$\{[(b_1+e_x) : (d_1+e_z)], e_x, e_z\in \{0,1\}\}\cap\{[(b_2+e_x) : (d_2+e_z)], e_x, e_z\in \{0,1\}\}=\emptyset.$$
\end{itemize}
Si $k=\mathbb Q$ ou $k=\FF_p$ un corps fini avec $p\geq 13$, on peut prendre par exemple $$l_1=x+y+2z,\, l_2=3x+3y+z.$$}

\medskip

\prop\label{ad}{Soit $F=\FF_p(x,y)$ le corps des fractions rationnelles \`a deux variables sur le corps fini $\FF_p$. Soit  $Q$ la quadrique lisse dans $\mathbb P^4_{F}$ d'\'equation homog\`ene
$$x_0^2-ax_1^2-fx_2^2+afx_3^2-g_1g_2x_4^2=0$$ avec $a\in \FF_p^*\setminus \FF_p^{*2}$ et $f, g_1, g_2$ d\'efinis comme dans (\ref{coef}) pour $k=\FF_p$. Alors le groupe $H^3_{\mathrm{nr}}(Q,\Z)$ est non nul.}
\proof{
Notons $A_x$, resp. $A_y$, resp. $A_z$, resp. $B_{i,j}$, l'anneau de valuation discr\`ete associ\'e au point g\'en\'erique de la droite $x=0$, resp. $y=0$ resp. $z=0$, resp. $l_i+h_j=0$,  $i=1,2$, $j=1,\ldots, 8$.

Il s'agit de v\'erifier les conditions $1$, $2$ et $3$ de la proposition \ref{cond}. Soit $B$ un anneau de valuation discr\`ete de $F$. On a les cas suivants \`a consid\'erer :

\begin{enumerate}
\item $B$ correspond \`a un point de codimension $1$ de  $\PF^2$.
\begin{enumerate}

\item Si $B$ est diff\'erent de $A_x, A_y, A_z, B_{i,j}$, le r\'esidu  $\partial_B(a, f, g_i)$, $i=1,2$, est nul, puisque les fonctions $a, f, g_1, g_2$ sont inversibles dans un tel anneau $B$.

\item $B=B_{i,j}$. Si $r\neq i$, $r=1,2$, alors $\partial_{B_{r,j}}(a, f, g_i)=0$ comme le cas pr\'ec\'edent. Fixons $i\in\{1,2\}$. Montrons que $\partial_{B_{i,j}}(a, f, g_i)\neq 0$. Supposons $h_j=0$, les autres cas sont identiques. Soit $k$ le corps r\'esiduel de $B_{i,0}$, i.e. le corps des fonctions de la droite $b_ix+c_iy+d_iz=0$, $b_i, c_i,d_i\in k^*$ (pour $h_j$ diff\'erent de z\'ero on utilise ainsi l'hypoth\`ese que $b_i, c_i, d_i$ sont diff\'erents de $-1$).  D'apr\`es le lemme \ref{prres}.1, $\partial_{B_{i,0}}(a, f, g_i)= (a, \frac xy)\in H^2(k,\Z)$. Apr\`es  passage \`a des coordonn\'ees affines, on est r\'eduit \`a \'etablir que le cup-produit $(a,x)$ n'est pas nul dans $H^2(\FF_p(x),\Z)$. On le voit par exemple en appliquant le lemme \ref{prres}.1 \`a l'anneau de valuation discr\`ete associ\'e \`a $x=0$ : $a\in \FF_p$ est non carr\'e.

\item $B=A_x$. Montrons que $\partial_{A_x}(a, f, g_i)=0$, $i=1,2$. Soit $k_x$ le corps r\'esiduel de $A_x$, i.e. le corps des fonctions de la droite  $x=0$. D'apr\`es le lemme \ref{prres}.1, $\partial_{A_x}(a, f, g_i)=-\partial_{A_x}(a, g_i, f)=-(a, g_{i,x})\in H^2(k_x, \Z))$, où $g_{i,x}$  d\'esigne la fonction induite par $g_i$ sur la droite $x=0$. Mais  $g_{i,x}$ est un carr\'e dans $k_x$, d'où $\partial_{A_x}(a, f, g_i)=0$ d'apr\`es \ref{prres}.2.

\item $B=A_y$. Comme dans le cas pr\'ec\'edent,  $\partial_{A_y}(a, f, g_2)=(a,g_{2,y})=0$.

\item $B=A_z$. Alors $\partial_{A_z}(a, f, g_1)=0$, car les fonctions $a, f, g_1$ sont inversibles dans $A_z$.
\end{enumerate}

\item $B$ correspond \`a un point ferm\'e $M$ de  $\PF^2$.
\begin{enumerate}
\item Si $M$ n'est pas situ\'e sur une des deux droites $x=0$, $y=0$, alors $f$ est inversible dans $B$.
\item Si $M$ est situ\'e sur une des deux droites $x=0$, $y=0$, alors l'une au moins des fonctions $g_1\frac{y^8}{z^8}$, $g_2$ est inversible dans $B$ d'apr\`es les hypoth\`eses (ii)-(iii), car le syst\`eme $xy=0, \prod_j(l_1+h_j)=0, \prod_j(l_2+h_j)=0$ n'a pas de solutions.
\end{enumerate}
\end{enumerate}
\qed
}

On finit par d\'ecrire explicitement  les exemples \'enonc\'es.

\theo\label{main}{Soit  $Q$ la quadrique lisse dans $\mathbb P^4_{\mathbb Q(x,y)}$ d'\'equation homog\`ene
$$x_0^2-ax_1^2-fx_2^2+afx_3^2-g_1g_2x_4^2=0$$ avec $a\in \mathbb Q^*\setminus \mathbb Q^{*2}$ et $f, g_1, g_2$ d\'efinis comme dans (\ref{coef}) pour $k=\mathbb Q$. Soit $X$ un mod\`ele projectif et lisse de $Q$ sur  $\mathbb P^2_{\mathbb Q}$ : $X$ est une $k$-vari\'et\'e projective et lisse et admet une fibration sur $\mathbb P^2_{\mathbb Q}$ \`a fibre g\'en\'erique $Q$.  Pour  une infinit\'e de nombres premiers $p$, la r\'eduction $X_p$ de $X$ modulo $p$ est bien d\'efinie et est une $\FF_p$-vari\'et\'e projective et lisse, telle que :
\begin{itemize}
\item[(i)] $\Hn(X_p, \Z)\neq 0$;
\item[(ii)] l'application $CH^2(X)\to CH^2({\bar X})^G$ n'est pas surjective.
\end{itemize}}
\proof{ D'apr\`es Hironaka, un mod\`ele projectif et lisse $X$ de $Q$ comme dans l'\'enonc\'e existe. De plus, pour une infinit\'e de nombres premiers $p$, l'image de $a$ dans $\FF_p$ n'est pas un carr\'e (par Chebotarev, ou par application de la loi de r\'eciprocit\'e quadratique) et la vari\'et\'e $X$ a  bonne r\'eduction en $p$ : $X_p$ est lisse. D'apr\`es la proposition \ref{ad}, le groupe $H^3_{\mathrm{nr}}(X_p,\Z)$ est non nul. Le th\'eor\`eme \ref{comp} permet de conclure. \qed\\}


\begin{thebibliography}{25}
\bibitem[A]{A} J.Kr. Arason,
\emph{Cohomologische invarianten quadratischer Formen},
J. Algebra \textbf{36} (1975), no. 3, 448--491.
\bibitem[CT]{CT} J.-L. Colliot-Th\'el\`ene, \emph{Birational invariants, purity and the Gersten conjecture},  $K$-theory and algebraic geometry: connections with quadratic forms and division algebras (Santa Barbara, CA, 1992),  1--64,
Proc. Sympos. Pure Math., \textbf{58}, Part 1, Amer. Math. Soc., Providence, RI, 1995.
\bibitem[CTOj]{CTOj} J.-L. Colliot-Th\'el\`ene et M. Ojanguren, \emph{Vari\'et\'es unirationnelles non rationnelles: au-del\`a de l'exemple d'Artin et Mumford},
Invent. Math. \textbf{97} (1989), no. 1, 141--158.

\bibitem[CTR]{CTR} J.-L. Colliot-Th\'el\`ene et W. Raskind, \emph{$K_2$-cohomology and the second Chow group},
Math. Ann. \textbf{270} (1985), no. 2, 165--199.
\bibitem[CTV]{CTV} J.-L. Colliot-Th\'el\`ene et C. Voisin, \emph{Cohomologie non ramifi\'ee et conjecture de Hodge enti\`ere}, en pr\'eparation.
\bibitem[G]{G} T. Geisser, \emph{Bass's conjectures and Tate's conjecture over finite fields}, en pr\'eparation.
\bibitem[K]{K} B. Kahn, \emph{Applications of weight-two motivic cohomology},
Doc. Math. \textbf{1} (1996), No. 17, 395--416.
\bibitem[KS]{KS} K. Kato and S. Saito, \emph{Unramified class field theory of arithmetical surfaces},
Ann. of Math. (2) \textbf{118} (1983), no. 2, 241--275.
\bibitem[LW]{LW} S. Lang et A. Weil,
\emph{Number of points of varieties in finite fields},
Amer. J. Math. \textbf{76}, (1954). 819--827.
\bibitem[M]{M} J.S. Milne, \emph{Zero cycles on algebraic varieties in nonzero characteristic: Rojtman's theorem},
Compositio Math. \textbf{47} (1982), no. 3, 271--287.
\bibitem[Q]{Q} D. Quillen, \emph{Higher algebraic $K$-theory I}, Algebraic $K$-theory, I: Higher $K$-theories (Proc. Conf., Battelle Memorial Inst., Seattle, Wash., 1972), pp. 85--147. Lecture Notes in Math., Vol. \textbf{341}, Springer, Berlin 1973.
\bibitem[R]{R}  A. A. Rojtman, \emph{The torsion of the group of $0$-cycles modulo rational equivalence},
Ann. of Math. (2) \textbf{111} (1980), no. 3, 553--569.
\bibitem[S]{Se} J-P. Serre, \emph{Corps locaux},  Publications de l'Universit\'e de Nancago, No. VIII. Hermann, Paris, 1968.
\bibitem[Sh]{Sh} C. Sherman,
\emph{$K$-cohomology of regular schemes},
Comm. Algebra \textbf{7} (1979), no. 10, 999--1027.
\end{thebibliography}
\end{document}